\documentclass[12pt,reqno]{amsart}
\usepackage{amsfonts,amsmath,amssymb}

\allowdisplaybreaks
\advance\textheight by \topskip

\usepackage{mathrsfs}

\usepackage[latin1]{inputenc}
\usepackage{graphicx}

\usepackage{tikz}

\usepackage{euscript}

\def\[{[\! [}
\def\]{]\! ]}

\begin{document}

\title{Partial sums of Horadam sequences: sum-free representations via generating functions}

\author[H.~Prodinger]{Helmut Prodinger}

\address{Helmut Prodinger,
Mathematics Department, Stellenbosch University,
7602 Stellenbosch, South Africa.}
\email{hproding@sun.ac.za}

\date{\today}

\begin{abstract}
Horadam sequences and their partial sums are computed via generating functions. The results are as simple as possible.
\end{abstract}

\maketitle

\section{Introduction}

Horadam sequences \cite{Horadam, Cooper} $W_n=W_n(a,b;p,q)$ are defined via
\begin{equation*}
W_n=pW_{n-1}+qW_{n-2}, \quad n\ge2,\quad W_0=a,\ W_1=b.
\end{equation*}
These numbers are of course a generalization of Fibonacci numbers, Lucas numbers and many others.
The characteristic equation 
\begin{equation*}
X^2-pX -q=0
\end{equation*}
is  essential, and the two roots are
\begin{equation*}
\lambda=\frac{p+\sqrt{p^2+4q}}{2},\quad \mu=\frac{p-\sqrt{p^2+4q}}{2}.
\end{equation*}
We define
\begin{equation*}
F_n=\frac{\lambda^n-\mu^n}{\lambda-\mu}\quad\text{and}\quad L_n=\lambda^n+\mu^n,
\end{equation*}
as these sequences resemble Fibonacci resp.\ Lucas numbers, and each solution of the recursion may be expressed
as a linear combination of these two sequences. We have $F_0=0$, $F_1=1$, $L_0=2$, $L_1=p$. For instance
\begin{equation*}
W_n=\bigl(b-\frac{ap}{2}\bigr)F_n+\frac a2 L_n.
\end{equation*}
The paper \cite{Cooper} concentrates on finding expressions for 
\begin{equation*}
\sum_{n\le k\le n+m}W_k=\sum_{0\le k\le n+m}W_k-\sum_{0\le k\le n-1}W_k.
\end{equation*}
In the rest of this short paper, we will find simple expressions for 
\begin{equation*}
S_n:=\sum_{0\le k\le n}W_k
\end{equation*}
using generating functions. The results do not contain summations, and can be expressed with the sequences $F_n$ and
$L_n$.

\section{Generating functions}
Standard computations produce
\begin{equation*}
W(z)=\sum_{k\ge0}W_kz^k=\frac{a+z(b-pa)}{1-pz-qz^2};
\end{equation*}
furthermore
\begin{equation*}
	F(z)=\sum_{k\ge0}F_kz^k=\frac{z}{1-pz-qz^2} \quad\text{and}\quad
L(z)=\sum_{k\ge0}L_kz^k=\frac{2-pz}{1-pz-qz^2}.
\end{equation*}
By general principles,
\begin{align*}
S(z)&=\sum_{k\ge0}S_kz^k=\frac{1}{1-z}\frac{a+z(b-pa)}{1-pz-qz^2}\\
&=\frac{a-pa+b}{1-p-q}\frac1{1-z}+\frac{-b-qa+qz(pa-a-b)}{1-p-q}\frac{1}{1-pz-qz^2}\\
&=\frac{a-pa+b}{1-p-q}\frac1{1-z}\\&
-\frac{2qa-pqa+2qb+pb}{2(1-p-q)}\frac{z}{1-pz-qz^2}
-\frac{qa+b}{2(1-p-q)}\frac{2-pz}{1-pz-qz^2}.
\end{align*}
Reading off the coefficient of $z^n$ on both sides leads to
\begin{align*}
S_n=\frac{a-pa+b}{1-p-q}-\frac{2qa-pqa+2qb+pb}{2(1-p-q)}F_n-\frac{qa+b}{2(1-p-q)}L_n.
\end{align*}
This answers the finite sum problem addressed in \cite{Cooper} completely, since the answer is
\begin{equation*}
S_{n+m}-S_{n-1}=-\frac{2qa-pqa+2qb+pb}{2(1-p-q)}(F_{n+m}-F_{n-1})
-\frac{qa+b}{2(1-p-q)}(L_{n+m}-L_{n-1}).
\end{equation*}

\bibliographystyle{plain}

\end{document}